\documentclass[12pt]{article}
\usepackage{amsthm,amsfonts, amsbsy, amssymb,amsmath,graphicx}
\usepackage{graphics}

\newtheorem{theorem}{Theorem}[section]
\newtheorem{lemma}{Lemma}[section]
\newtheorem{corollary}{Corollary}[section]

\newtheorem{statement}{Statement}[section]

 \def\thtext#1{
 \catcode`@=11
 \gdef\@thmcountersep{. #1}
 \catcode`@=12}

 \newcounter{Rk}[section]
 \renewcommand{\thtext}{\thesection.\arabic{Rk}}
 \newenvironment{remark}{\trivlist\item[\hskip\labelsep{\bf Remark}]
 \refstepcounter{Rk}{\bf\thesection.\arabic{Rk}.}}%
 {\endtrivlist}

 \newcounter{Df}[section]
 \renewcommand{\thtext}{\thesection.\arabic{Df}}
 \newenvironment{definition}{\trivlist\item[\hskip\labelsep{\bf Definition}]
 \par\refstepcounter{Df}{\bf\thesection.\arabic{Df}.}}%
 {\endtrivlist}

 \newenvironment{example}{\trivlist \item[\hskip\labelsep{\bf Example.}]}%
 {\endtrivlist}

\title{Graph-Links}

\author{Denis Petrovich Ilyutko\footnote{Partially supported by grants of RF President
NSh -- 660.2008.1, RFBR 07--01--00648, RNP 2.1.1.3704, the Federal
Agency for Education NK-421P/108.}, Vassily Olegovich
Manturov\footnote{Partially supported by grants of RFBR
07--01--00648}}

\begin{document}

\maketitle

\abstract{The present paper is a review of the current state of {\em
Graph-Link Theory} (graph-links are also closely related to {\em
homotopy classes of looped interlacement graphs}): theory suggested
in~\cite{IM1,IM2}, see also~\cite{TZ}, dealing with a generalisation
of knots obtained by translating the Reidemeister moves for links
into the language of intersection graphs of chord diagrams. In this
paper we show how some methods of classical and virtual knot theory
can be translated into the language of abstract graphs, and some
theorems can be reproved and generalised to this graphical setting.
We construct various invariants, prove certain minimality theorems
and construct functorial mappings for graph-knots and graph-links.
In this paper, we first show non-equivalence of some graph-links to
virtual links.

}

\section{Introduction}

\begin{figure}
\centering\includegraphics[width=200pt]{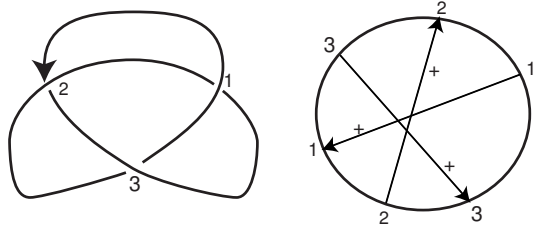} \caption{The
Right Trefoil and Its Gauss Diagram}\label{Gaussdiagram}
 \end{figure}

\begin{figure}
\centering\includegraphics[width=200pt]{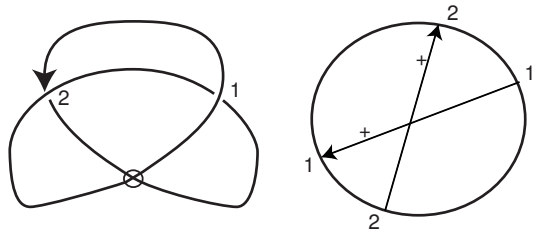} \caption{The
Virtual Trefoil and Its Gauss Diagram} \label{VGaussdiagram}
 \end{figure}

It is well known that classical and virtual knots can be represented
by Gauss diagrams, and the whole information about the knot and its
invariants can be read out of any Gauss diagram encoding it, see
Fig.~\ref{Gaussdiagram}. Whenever a Gauss diagram does not describe
any {\em embedded curve} in $\mathbb{R}^{2}$ (just because the
corresponding Gauss code is not planar) as in
Fig.~\ref{VGaussdiagram}, one gets a virtual knot; the generic
immersion point is encircled.

\begin{figure}
  \centering\includegraphics[width=200pt]{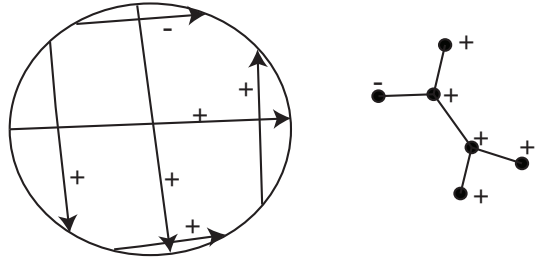}
  \caption{A Gauss Diagram and its Labeled Intersection Graph}
  \label{gaussdiagram}
 \end{figure}

\begin{figure}
\centering\includegraphics[width=200pt]{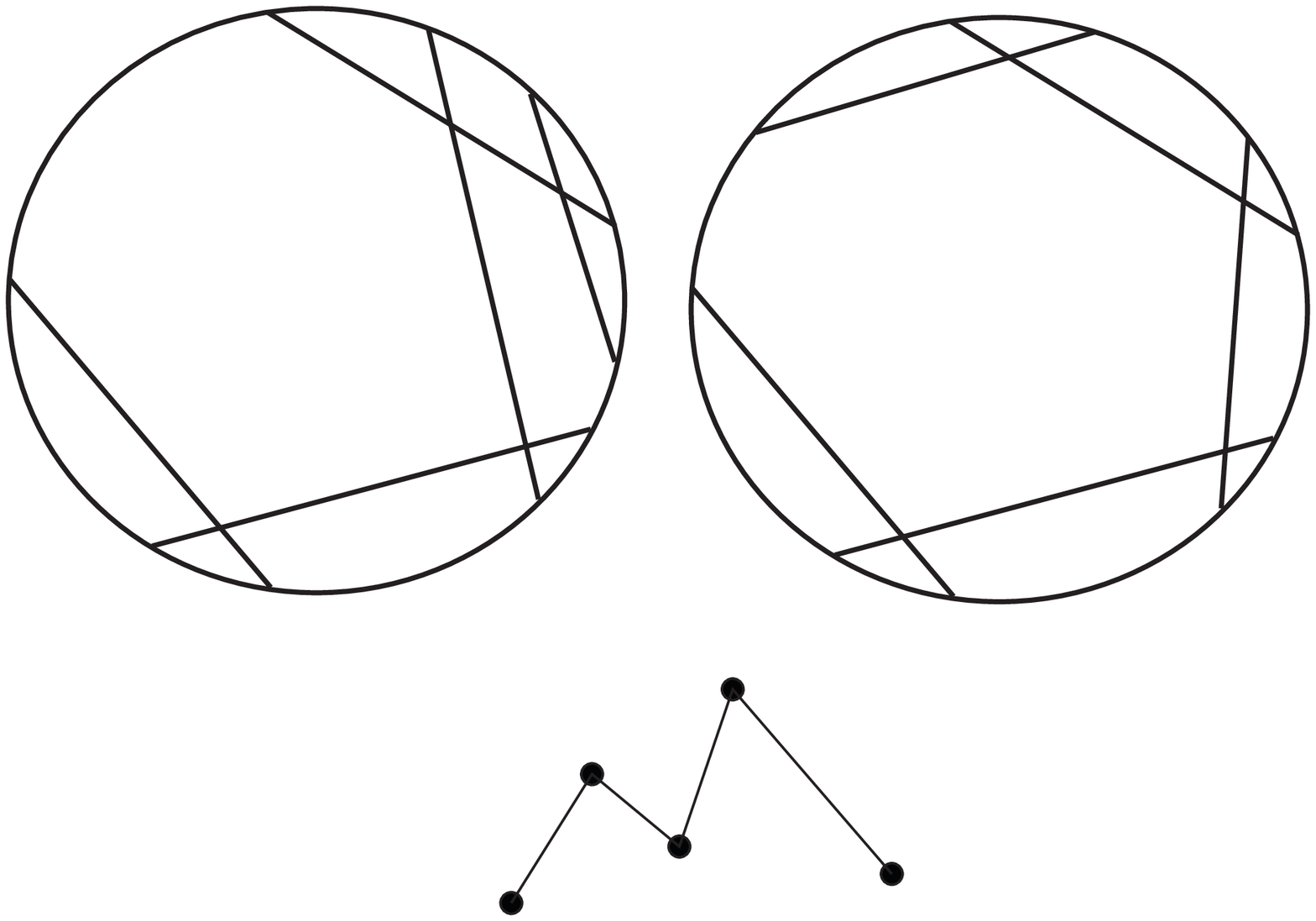}\caption{A
Graph not Uniquely Represented by Chord Diagrams}
\label{nonuniqueness}
 \end{figure}

\begin{figure}
\centering\includegraphics[width=300pt]{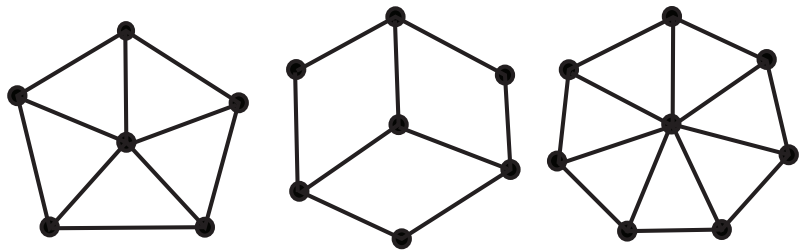}
\caption{Non-Realisable Bouchet Graphs} \label{Bouchet}
\end{figure}

It turns out that some information about the knot can be obtained
from a more combinatorial data: intersection graphs of Gauss
diagrams. The intersection (adjacency) graph is a {\em simple
graph}, i.e.\ a graph without loops and multiple edges, whose
vertices are in one-to-one correspondence with chords of the Gauss
diagram of the knot (the latter are, in turn, in one-to-one
correspondence with crossings of the knot). Two vertices of the
intersection graph are adjacent whenever the corresponding arrows of
the Gauss diagram are {\em linked}, i.e.\ the corresponding chords
``intersect'' each other in the picture, see
Fig.~\ref{gaussdiagram}. Vertices of the intersection graph are
endowed with the local writhe number of the crossing. However,
sometimes the Gauss diagram can be obtained from the intersection
graph in a non-unique way, see Fig.~\ref{nonuniqueness}, and some
graphs can not be represented by chord diagrams at all~\cite {Bou},
see Fig.~\ref{Bouchet}.

When passing to the intersection graph, we remember the writhe
number information, but forget about the {\em right-left}
information encoded by the arrows. Principally, it is possible to
describe analogous objects when all information is saved in the
intersection graph; however, already the writhe number information
is sufficient to recover a lot of data, as we shall see. Even more,
if we forget about the writhe number information and only have the
structure of opposite edges we shall get non-trivial objects (modulo
Reidemeister's moves).

 \begin{figure}
\centering\includegraphics[width=200pt]{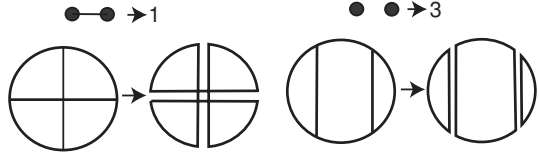}
\caption{Resmoothing along Two Chords Yields One or Three Circles}
\label{resmoothings}
 \end{figure}

Probably, the simplest evidence that one can get some information
out of the intersection graph is the number of circles one gets in a
certain state after a smoothing, see Fig.~\ref{resmoothings}. The
simplest example shows that the number of circles obtained after
resmoothings in two crossings gives three circles if the
corresponding chords are unlinked or one circle if they aren't. The
general theorem (Soboleva's theorem, see ahead) allows one to count
the number of circles in Kauffman's states out of the intersection
graph. In particular, this means that graphs not necessarily
corresponding to any knot admit a way of generalising the Kauffman
bracket, which coincides with the usual Kauffman bracket when the
graph is realisable by a knot. This was the initial point of
investigation for L.~Traldi and L.~Zulli~\cite{TZ} (looped
interlacement graphs): they constructed a self-contained theory of
``non-realisable knots'' possessing lots of interesting knot
theoretic properties. These objects are equivalent classes of
(decorated) graphs modulo ``Reidemeister moves'' (translated into
the language of intersection graphs). A significant disadvantage of
this approach was that it had applications only to knots, {\em not
links}: in order to encode a link, one has to use a more complicated
object rather than just a Gauss diagram, a Gauss diagram {\em on
many circles}. This approach was further developed in Traldi's
works~\cite{Tr1,Tr2}, and it allowed to encode not only knots but
also links with any number of components by decorated graphs. The
important question arises here: whether or not every graph is
Reidemeister equivalent to the looped interlacement graph of a
virtual knot diagram.

\begin{figure}
\centering\includegraphics[width=300pt]{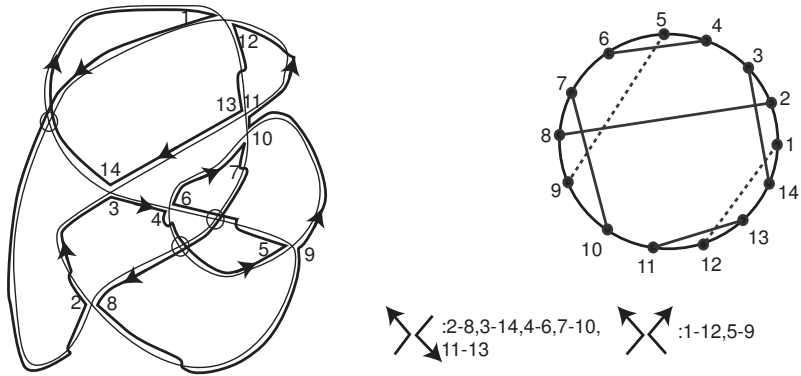} \caption{Rotating
Circuit Shown by a Thick Line; Chord Diagram} \label{rotations}
 \end{figure}

We suggested another way of looking at knots and links and
generalising them: whence a Gauss diagram corresponds to a
transverse passage along a knot, one may consider a {\em rotating
circuit} which never goes straight and always turns right or left at
a classical crossing. One can also encode the type of smoothing
(Kauffman's $A$-smoothing or Kauffman's $B$-smoothing) corresponding
to the crossing where the circuit turns right or left and never goes
straight, see Fig.~\ref{rotations}. We note that chords of the
diagrams are naturally split into two sets: those corresponding to
crossings where two opposite directions correspond to {\em
emanating} edges with respect to the circuit and the other two
correspond to {\em incoming} edges, and those where we have two
consecutive (opposite) edges one of which is incoming and the other
one is emanating.

 \begin{figure}
\centering\includegraphics[width=200pt]{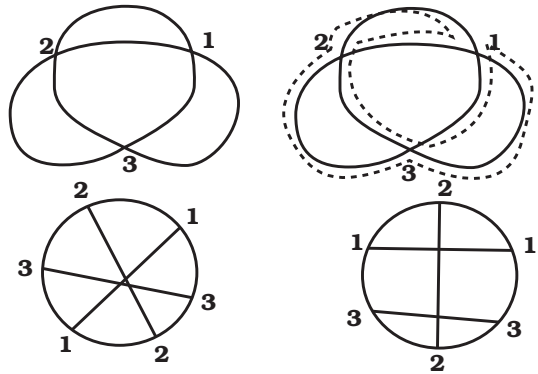} \caption{A
Gauss Circuit; A Rotating Circuit} \label{Gaussrot}
 \end{figure}

Certainly, there is an intuitive way of transforming Gauss diagrams
into rotating circuits and vice versa (in the case of knots), see
Fig.~\ref{rotations},~\ref {Gaussrot}. The second important question
arises here whether there is an equivalence between the set of
homotopy classes of looped interlacement graphs introduced by
L.~Traldi and L.~Zulli~\cite{TZ} and the set of graph-knots
introduced by ourselves~\cite {IM1,IM2}. The equivalence of these
two sets was proved in~\cite {Ily1,Ily2} and, moreover, the homotopy
class of looped interlacement graphs and the graph-knot both
constructed from a given virtual knot diagram are related by this
equivalence. This construction will be described in detail in
Section~\ref {gllig}. Moreover, one can consider mixed circuits
(which was initiated by Traldi~\cite{Tr1}). In this paper we do not
consider mixed circuits. The way of translating ``rotating
circuits'' into ``transverse circuits'' can be rewritten in a
combinatorial way and translated into the language of intersection
graphs. In some sense it is obvious that {\em whenever an
intersection graph is realisable in the sense of Gauss diagrams, the
corresponding ``rotating'' intersection graph is realisable in the
sense of rotating diagrams}, just because if one of these two graphs
is realisable, the corresponding 4-valent graph can be just drawn on
the plane (with virtual crossings) and the algorithm becomes a
``real redrawing algorithm''. This statement follows immediately
from our equivalence. This allows us to switch from Gauss diagrams
to rotating diagrams and vice versa, whenever we are proving some
non-realisability theorems. We shall answer the question of whether
or not every graph is Reidemeister equivalent to the looped
interlacement graph of a virtual knot diagram, and by using the
equivalence we immediately get the answer to this question
concerning graph-knots. We shall show that it is not true for looped
interlacement graphs and by using the equivalence it is not true for
graph-links. The examples of ``non-realisable'' looped interlacement
graphs appeared firstly in the papers~\cite {Ma1,Ma2}. So, the
theory of graph-links is interesting for various reasons:

a) In some cases it exhibits purely combinatorial ways of extracting
invariants for knots.

b) In some cases it produces heuristic approaches to new ``knot
theories''.

c) It highlights some ``graphical'' effects which are hardly visible
in usual or virtual knot theory.

\begin{figure}
\centering\includegraphics[width=100pt]{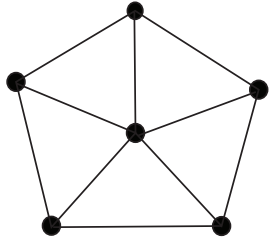}
\caption{The First Bouchet Graph Gives a Non-Realisable Graph-Knot}
\label{Bouchet9}
 \end{figure}

We conclude the introduction part by a couple of examples of free
knots which are not equivalent to any realisable knot, see
Fig.~\ref{Bouchet9}. We call a graph-link {\em non-realisable} if it
has no realisable representative.

Here we do not indicate any crossing decoration because {\em any
graph-link with this underlying graph is non-realisable}. More
precisely, this example gives us a non-realisable {\em free
graph-link}. Having a virtual link, we may forget about over/under
and right/left information and take care only about the 4-valent
underlying graph with the structure of opposite edges. This leads to
the notion of {\em free knots} and {\em free links} considered
in~\cite {Ma1, Ma2}. It turns out that some information about
virtual knots can be caught just from the 4-valent graph, which
proves non-triviality of many free knots and free links. The same
trick works for graph-links: however, here instead of the underlying
four-valent graph we consider the abstract graph which plays the
role of the intersection graph of the {\em non-existing chord
diagram}.

\begin{figure}
\centering\includegraphics[width=250pt]{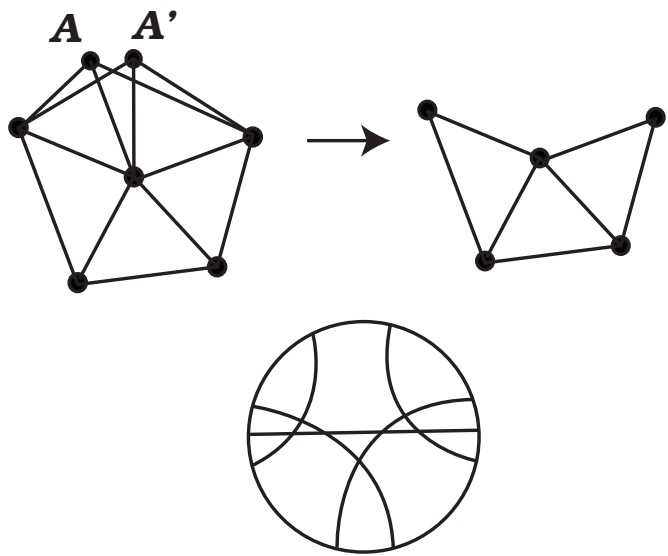} \caption{A
Non-Realisable Graph Representing a Trivial Graph-Link}
\label{Bouchet10}
 \end{figure}

Certainly, the graph shown in Fig.~\ref{Bouchet10} (left upper) is
itself non-realisable (in the sense of looped interlacement graphs
and Gauss diagrams), but what if we decorate its crossings in some
way and then try to apply Reidemeister moves hoping to make it
realisable. For some graph-links it is possible, see, e.g.\
Fig.~\ref{Bouchet10}.

The graph-link with ``Gauss diagram'' shown in Fig.~\ref{Bouchet10}
(left upper) is realisable. Indeed, the second Reidemeister move
translated into the language of Gauss diagram is an addition/removal
of two ``parallel'' chords. In the language of intersection graphs,
chords correspond to vertices, and ``parallel'' chords correspond to
vertices having the same set of adjacent vertices. So, the vertices
$A$ and $A'$ in Fig.~\ref{Bouchet10} are adjacent, and the removal
of these two vertices makes our diagram realisable.

\begin{figure}
\centering\includegraphics[width=200pt]{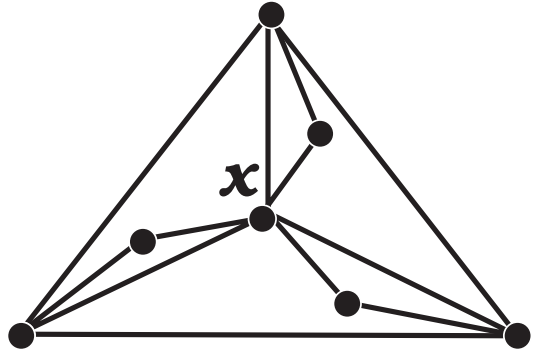} \caption{An Even
Non-Realisable Graph Representing a Non-Realisable Graph-Knot}
\label{xx}
 \end{figure}

The problem of finding free links having no representative
realisable by a chord diagram is a problem similar to the problem of
constructing virtual knots not equivalent to any classical link.
Surprisingly, the solution to the problem in the case of free links
can be achieved by using {\em parity considerations} (Theorem~\ref
{th:main_pair}): all the crossings shown in Fig.~\ref{Bouchet9} are
{\em odd} (each of them is adjacent to an odd number of other
crossings) and there is no immediate way to contract any two them by
using a second Reidemeister moves. This is indeed sufficient for a
graph-link to be non-realisable in a very strong sense: any diagram
of this link has {\em a spur} of the initial diagram (we disregard
the writhe number information), which is non-realisable, and, the
graph-link is, in turn, itself non-realisable. Also we have an
example of non-realisable graph all the vertices of it are even, see
Fig.~\ref {xx}.

In fact, parity arguments allow to prove very strong theorems in the
realm of virtual knots, and we strongly recommend the reader to read
the paper~\cite{Ma3} in the present volume.

The present paper is organized as follows.

We first give definitions of graph-links from two points of view:
rotating circuits and Gauss diagrams, and describe their
interactions.

In the third section we introduce parity and prove non-triviality
results. In particular, parity arguments allow one to construct
graph-valued invariants of graph-links.

The fourth section is devoted to the {\em orientability} aspects of
$2$-surfaces corresponding to virtual knots (atoms) and functorial
mappings for virtual knots and free links defined by using parities.

In section 5, we briefly describe the way of extending the Kauffman
bracket and some other invariants. We also formulate some minimality
results for graph-links and conclude the paper by a list of unsolved
problems.

\section*{Acknowledgments}

The authors are grateful to L.\,H.~Kauffman, V.\,A.~Vassiliev, and
A.\,T.~Fomenko for their interest to this work.

\section{Graph-Links and Looped Interlacement Graphs}\label {gllig}

\subsection{Chord diagrams and Framed 4-Graphs}

Throughout the paper all graphs are finite. Let $G$ be a graph with
the set of vertices $V(G)$ and the set of edges $E(G)$. We think of
an edge as an equivalence class of the two half-edges. We say that a
vertex $v\in V(G)$ has {\em degree} $k$ if $v$ is incident to $k$
half-edges. A graph whose vertices have the same degree $k$ is
called {\em $k$-valent} or a {\em k-graph}. The free loop, i.e.\ the
graph without vertices, is also considered as $k$-graph for any $k$.

 \begin {definition}\label {def:fr4}
A $4$-graph is {\em framed} if for every vertex the four emanating
half-edges are split into two pairs of (formally) opposite edges.
The edges from one pair are called {\em opposite to each other}.
 \end {definition}

\begin {figure}
\centering\includegraphics[width=250pt]{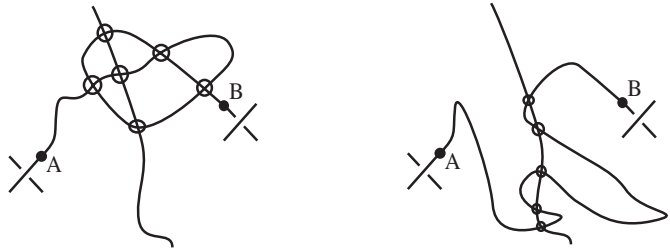} \caption{The
detour move} \label{detour}
 \end {figure}

A {\em virtual diagram} is a framed $4$-graph embedded into
${\mathbb R}^2$ where each crossing is either endowed with a
classical crossing structure (with a choice for underpass and
overpass specified) or just said to be virtual and marked by a
circle. A {\em virtual link} is an equivalence class of virtual
diagrams modulo generalised Reidemeister moves. The latter consist
of usual Reidemeister moves referring to classical crossings and the
{\em detour move} that replaces one arc containing only virtual
intersections and self-intersections by another arc of such sort in
any other place of the plane, see Fig.~\ref{detour}. A {\em
projection} of a virtual diagram is a framed $4$-graph obtained from
the diagram by considering classical crossings as vertices and
virtual crossings are just intersection points of images of
different edges. A virtual diagram is {\em connected} if its
projection is connected. Without loss of generality, {\it all
virtual diagrams are assumed to be connected and contain at least
one classical crossing}~\cite {IM1,IM2}.

  \begin{definition}
A {\em chord diagram} is a cubic graph consisting of a selected
cycle (the {\em circle}) and several non-oriented edges ({\em
chords}) connecting points on the circle in such a way that every
point on the circle is incident to at most one chord. A chord
diagram is {\em labeled} if every chord is endowed with a label
$(a,\alpha)$, where $a\in\{0,1\}$ is the framing of the chord, and
$\alpha\in\{\pm\}$ is the sign of the chord. If no labels are
indicated, we assume the chord diagram has all chords with label
$(0,+)$. Two chords of a chord diagram are called {\em linked} if
the ends of one chord lie in different connected components of the
circle with the end-points of the second chord removed.
  \end{definition}

\begin{figure}
\centering\includegraphics[width=200pt]{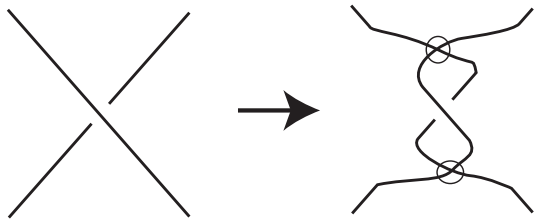}
\caption{Virtualisation} \label{virtualisation}
 \end{figure}

 \begin {definition}
By a {\em virtualisation} of a classical crossing of a virtual
diagram we mean a local transformation shown in
Fig.~\ref{virtualisation}.
 \end {definition}

Having a labeled chord diagram $D$, one can construct a virtual link
diagram $K(D)$ (up to virtualisation) as follows. Let us immerse
this diagram in ${\mathbb R}^{2}$ by taking an embedding of the
circle and placing some chords inside the circle and the other ones
outside the circle. After that we remove neighbourhoods of each of
the chord ends and replace them by a pair of lines (connecting four
points on the circle which are obtained after removing
neighbourhoods) with a classical crossing if the chord is framed by
$0$ and a couple of lines with a classical crossing and a virtual
crossing if the chord is framed by $1$ in the following way. The
choice for underpass and overpass is specified as follows. A
crossing can be smoothed in two ways: $A$ and $B$ as in the Kauffman
bracket polynomial; we require that the initial piece of the circle
corresponds to the $A$-smoothing if the chord is positive and to the
$B$-smoothing if it is negative:
 $A\colon\raisebox{-0.25\height}{\includegraphics[width=0.5cm]{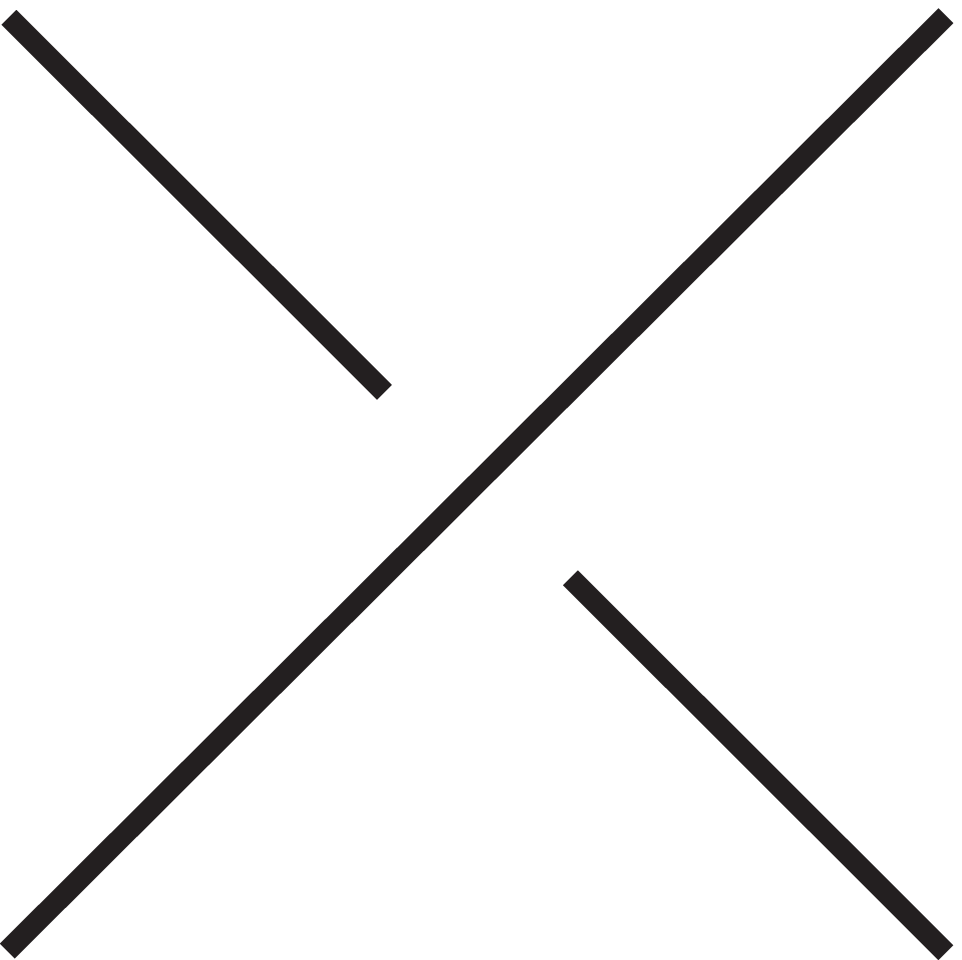}}\to
 \raisebox{-0.25\height}{\includegraphics[width=0.5cm]{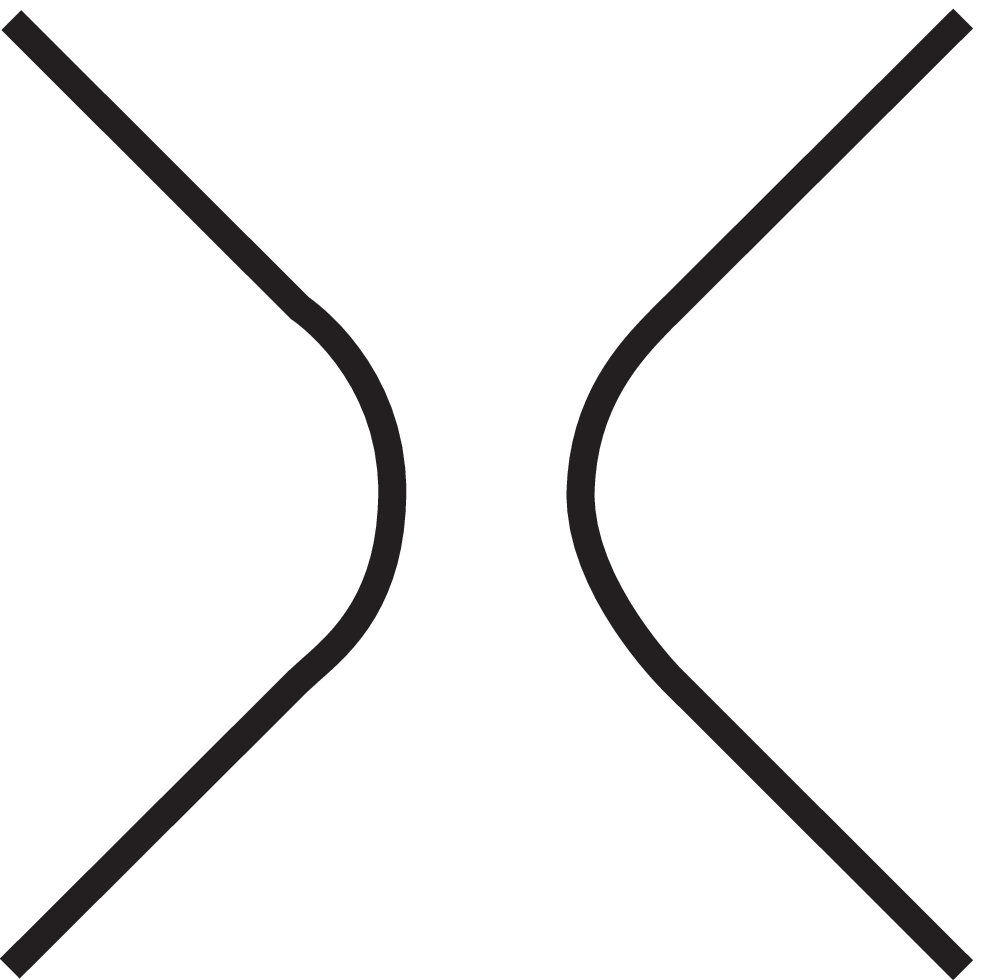}}$,
 $B\colon\raisebox{-0.25\height}{\includegraphics[width=0.5cm]{skcrossr.eps}}\to
 \raisebox{-0.25\height}{\includegraphics[width=0.5cm]{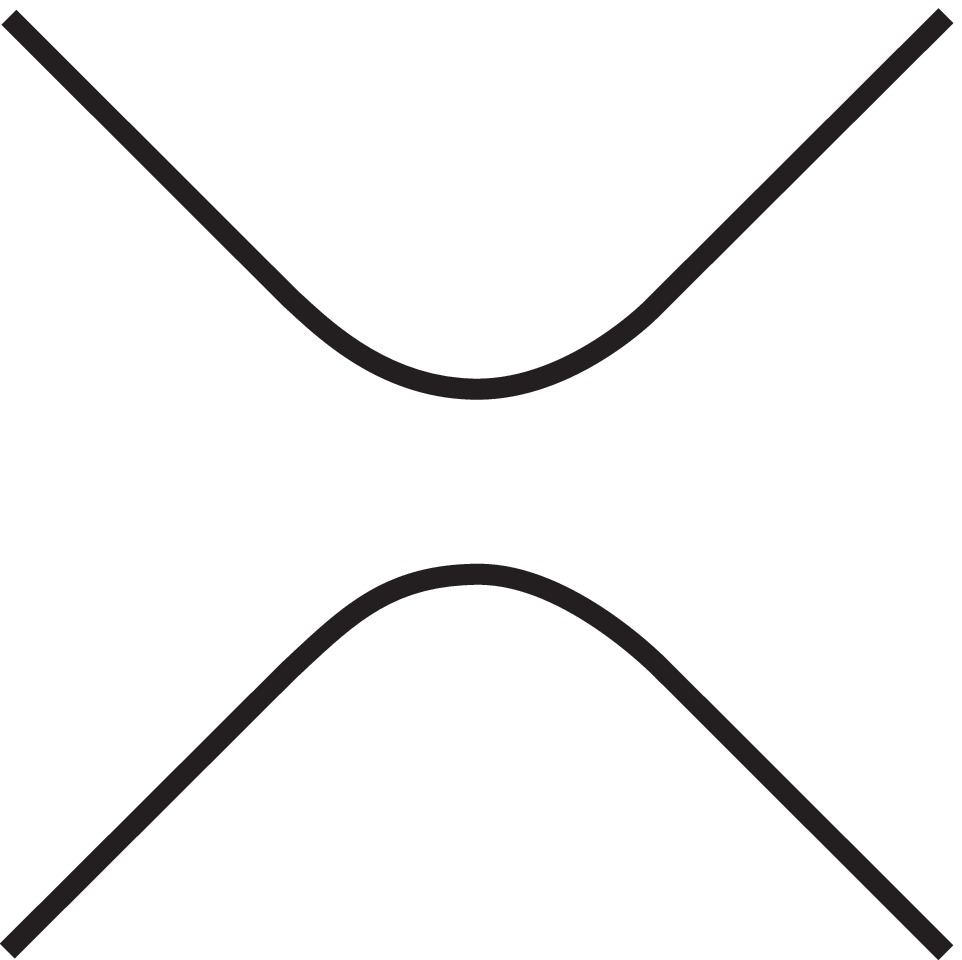}}$.

Conversely, having a connected virtual diagram $K$, one can get a
labeled chord diagram $D_C(K)$, see Fig.~\ref{rotations}. Indeed,
one takes a {\em circuit} $C$ of $K$ which is a map from $S^{1}$ to
the projection of $K$. This map is bijective outside classical and
virtual crossings, has exactly two preimages at each classical and
virtual crossing, goes transversally at each virtual crossing and
turns from an half-edge to an adjacent (non-opposite) half-edge at
each classical crossing. Connecting the two preimages of a classical
crossing by a chord we get a chord diagram, where the sign of the
chord is $+$ if the circuit locally agrees with the $A$-smoothing,
and $-$ if it agrees with the $B$-smoothing, and the framing of a
chord is $0$ (resp., $1$) if two opposite half-edges have the
opposite (resp., the same) orientation. It can be easily checked
that this operation is indeed inverse to the operation of
constructing a virtual link out of a chord diagram: if we take a
chord diagram $D$, and construct a virtual diagram $K(D)$ out of it,
then for some circuit $C$ the chord diagram $D_C(K(D))$ will
coincide with $D$. The rule for setting classical crossings here
agrees with the rule described above. This proves the following

  \begin {theorem}\cite {MBook}
For any connected virtual diagram $L$ there is a certain labeled
chord diagram $D$ such that $L=K(D)$.
  \end {theorem}

The Reidemeister moves on virtual diagrams generate the Reidemeister
moves on labeled chord diagrams~\cite {IM1, IM2}.

\subsection{Reidemeister Moves for Looped Interlacement Graphs and Graph-links}

Now we are describing moves on graphs obtained from virtual diagrams
by using rotating circuit~\cite {IM1,IM2} and the Gauss
circuit~\cite {TZ}. These moves in both cases will correspond to the
``real'' Reidemeister moves on diagrams. Then we shall extend these
moves to all graphs (not only to realisable ones). As a result we
get new objects, a {\em graph-link} and a {\em homotopy class of
looped interlacement graphs}, in a way similar to the generalisation
of classical knots to virtual knots: the passage from realisable
Gauss diagrams (classical knots) to arbitrary chord diagrams leads
to the concept of a virtual knot, and the passage from realisable
(by means of chord diagrams) graphs to arbitrary graphs leads to the
concept of two new objects, {\em graph-links} and {\em homotopy
classes of looped interlacement graphs} (here `looped' corresponds
to the writher number, if the writher number is -1 then the
corresponding vertex has a loop). To construct the first object we
shall use simple labeled graphs, and for the second one we shall use
(unlabeled) graphs without multiple edges, but loops are allowed.

 \begin {definition}
A graph is {\em labeled} if every vertex $v$ of it is endowed with a
pair $(a,\alpha)$,  where $a\in\{0,1\}$ is the framing of $v$, and
$\alpha\in\{\pm\}$ is the sign of $v$. Let $D$ be a labeled chord
diagram $D$. The {\em labeled intersection graph}, cf.~\cite{CDL},
$G(D)$ of $D$ is the labeled graph: 1) whose vertices are in
one-to-one correspondence with chords of $D$, 2) the label of each
vertex corresponding to a chord coincides with that of the chord,
and 3) two vertices are connected by an edge if and only if the
corresponding chords are linked.
 \end {definition}

 \begin {definition}
A simple graph $H$ is called {\em realisable} if there is a chord
diagram $D$ such that $H=G(D)$.
 \end {definition}

The following lemma is evident.

 \begin {lemma}\label {lem:co_re}
A simple graph is realisable if and only if each its connected
component is realisable.
 \end {lemma}

 \begin {definition}
Let $G$ be a graph and let $v\in V(G)$. The set of all vertices
adjacent to $v$ is called the {\em neighbourhood of a vertex} $v$
and denoted by $N(v)$ or $N_G(v)$.
 \end {definition}

Let us define two operations on simple unlabeled graphs.

 \begin {definition} (Local Complementation)
Let $G$ be a graph. The {\em local complementation} of $G$ at $v\in
V (G)$ is the operation which toggles adjacencies between $a,b\in
N(v)$, $a\ne b$, and doesn't change the rest of $G$. Denote the
graph obtained from $G$ by the local complementation at a vertex $v$
by $\operatorname{LC}(G;v)$.
 \end {definition}

 \begin {definition} (Pivot)
Let $G$ be a graph with distinct vertices $u$ and $v$. The {\em
pivoting operation} of a graph $G$ at $u$ and $v$ is the operation
which toggles adjacencies between $x,\,y$ such that
$x,\,y\notin\{u,v\}$, $x\in N(u),\,y\in N(v)$ and either $x\notin
N(v)$ or $y\notin N(u)$, and doesn't change the rest of $G$. Denote
the graph obtained from $G$ by the pivoting operation at vertices
$u$ and $v$ by $\operatorname{piv}(G;u,v)$.
 \end {definition}

\begin{figure}
\centering\includegraphics[width=200pt]{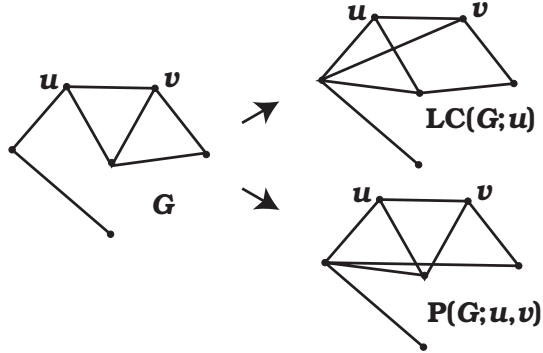} \caption{Local
Complementation and Pivot} \label{lcp}
 \end{figure}

 \begin {example}
In Fig.~\ref {lcp} the graphs $G$, $\operatorname{LC}(G;u)$ and
$\operatorname{piv}(G;u,v)$ are depicted.
 \end {example}

The following lemma can be easily checked.

 \begin {lemma}
If $u$ and $v$ are adjacent then there is an isomorphism
 $$
\operatorname{piv}(G;u,v)\cong\operatorname{LC}(\operatorname{LC}(\operatorname{LC}(G;u);v);u).
 $$
 \end {lemma}

Let us define graph-moves by considering intersection graphs of
chord diagrams constructed by using {\em a rotating circuit}, and
these moves correspond to the Reidemeister moves on virtual
diagrams. As a result we obtain a new object --- an equivalence
class of labeled graphs under formal moves. These moves were defined
in~\cite {IM1,IM2}.

 \begin {definition}
$\Omega_g 1$. The first Reidemeister graph-move is an
addition/removal of an isolated vertex labeled $(0,\alpha)$,
$\alpha\in\{\pm\}$.

$\Omega_g 2$. The second Reidemeister graph-move is an
addition/removal of two non-adjacent (resp., adjacent) vertices
having $(0,\pm\alpha)$ (resp., $(1,\pm\alpha)$) and the same
adjacencies with other vertices.

$\Omega_g 3$. The third Reidemeister graph-move is defined as
follows. Let $u,\,v,\,w$ be three vertices of $G$ all having label
$(0,-)$ so that $u$ is adjacent only to $v$ and $w$ in $G$. Then we
only change the adjacency of $u$ with the vertices $t\in
N(v)\setminus N(w)\bigcup N(w)\setminus N(v)$ (for other pairs of
vertices we do not change their adjacency). In addition, we switch
the signs of $v$ and $w$ to $+$. The inverse operation is also
called the third Reidemeister graph-move.

$\Omega_g 4$. The fourth graph-move for $G$ is defined as follows.
We take two adjacent vertices $u$ and $v$ labeled $(0,\alpha)$ and
$(0,\beta)$ respectively. Replace $G$ with
$\operatorname{piv}(G;u,v)$ and change signs of $u$ and $v$ so that
the sign of $u$ becomes $-\beta$ and the sign of $v$ becomes
$-\alpha$.

$\Omega_g 4'$. In this fourth graph-move we take a vertex $v$ with
the label $(1,\alpha)$. Replace $G$ with $\operatorname{LC}(G;v)$
and change the sign of $v$ and the framing for each $u\in N(v)$.
 \end {definition}

The comparison of the graph-moves with the Reidemeister moves yields
the following theorem.

 \begin {theorem}\label {th:link_graph}
Let $K_1$ and $K_2$ be two connected virtual diagrams, and let $G_1$
and $G_2$ be two labeled intersection graphs obtained from $K_1$ and
$K_2$, respectively. If $K_1$ and $K_2$ are equivalent in the class
of connected diagrams then $G_1$ and $G_2$ are obtained from one
another by a sequence of $\Omega_g 1 - \Omega_g 4'$ graph-moves.
 \end {theorem}

  \begin {definition}
A {\em graph-link} is an equivalence class of simple labeled graphs
modulo $\Omega_g 1 - \Omega_g 4'$ graph-moves.
 \end {definition}

  \begin {remark}
For a graph-link having representatives with orientable atoms, see
ahead, there are two formally different equivalence relations. The
first relation is described in~\cite {IM1} (which includes only
diagrams with orientable atoms) and the last one defines
graph-links. We do not know whether these equivalence relations
coincide for graph-links (they do for graph-knots, see ahead).
Nevertheless, we use the same term `graph-link' for the object
introduced in this paper.
 \end  {remark}

 \begin {remark}
Let us consider a simple realisable labeled graph. Let us represent
this graph as an intersection graph of a chord diagram $D$.
Constructing a virtual diagram $K(D)$ we have to restore the
structure of opposite edges at each vertex. The first components
(which are called the framings of vertices) of labels of vertices
are responsible for the structure of opposite edges. Therefore, we
shall consider framed graphs of two types. The first type is framed
4-graphs defined in Definition~\ref {def:fr4}. The second one is
{\em free framed graphs} which are equivalence classes of simple
labeled graphs with labels having only framings modulo $\Omega_g 4$
and $\Omega_g 4'$ graph-moves up to signs of labels (we disregard
the sign of each vertex).
 \end {remark}

 \begin {definition}
A {\em free graph-link} is an equivalence class of free framed
graphs modulo $\Omega_g 1 - \Omega_g 3$ graph-moves up to signs of
labels.
 \end {definition}

Let $D_G(K)$ be the Gauss diagram of a virtual diagram $K$. Let us
construct the graph obtained from the intersection graph of $D_G(K)$
by adding loops to vertices corresponding to chords with negative
writhe number~\cite {TZ}. We refer to this graph as the {\em looped
interlacement graph} or the {\em looped graph}. Let us construct the
moves on graphs. These moves are similar to the moves for
graph-links and also correspond to the Reidemeister moves on virtual
diagrams.

 \begin {definition}\label {def:mov_ga}
The first Reidemeister move for looped interlacement graphs is an
addition/removal of an isolated looped or unlooped vertex.

The second Reidemeister move for looped interlacement graphs is an
addition/removal of two vertices having the same adjacencies with
other vertices and, moreover, one of which is looped and the other
one is unlooped.

The third Reidemeister move for looped interlacement graphs is
defined as follows. Let $u,\,v,\,w$ be three vertices such that $v$
is looped, $w$ is unlooped, $v$ and $w$ are adjacent, $u$ is
adjacent to neither $v$ nor $w$, and every vertex
$x\notin\{u,\,v,\,w\}$ is adjacent to either $0$ or precisely two of
$u,\,v,\,w$. Then we only remove all three edges $uv$, $uw$ and
$vw$.  The inverse operation is also called the third Reidemeister
move.
 \end {definition}

The two third Reidemeister moves do not exhaust all the
possibilities for representing the third Reidemeister move on Gauss
diagrams~\cite {TZ}. It can be shown that all the other versions of
the third Reidemeister move are combinations of the second and third
Reidemeister moves, see~\cite {Ost} for details.

 \begin {definition}
We call an equivalence class of graphs (without multiple edges, but
loops are allowed) modulo the three moves listed in Definition~\ref
{def:mov_ga} a {\em homotopy class} of looped interlacement graphs.
A {\em free homotopy class} is an equivalence class of simple graphs
modulo the Reidemeister moves for looped interlacement graphs up to
loops, i.e.\ we forget about loops.
 \end {definition}

 \begin {remark}
Looped interlacement graphs encode only knot diagrams but
graph-links can encode virtual diagrams with any number of
components. The approach using a rotating circuit has an advantage
in this sense. In~\cite {Tr1} L.~Traldi introduced the notion of a
marked graph by considering any Euler tour (we have vertices which
we go transversally and in which we rotate).
 \end {remark}

\subsection{Looped Interlacement Graphs and Graph-Links}

Let $G$ be a labeled graph on vertices from the enumerated set
$V(G)=\{v_1,\dots,v_n\}$, and let $A(G)$ be the {\em adjacency
matrix} of $G$ over $\mathbb{Z}_2$ defined as follows: $a_{ii}$ is
equal to the framing of $v_i$, $a_{ij}=1$, $i\ne j$, if and only if
$v_i$ is adjacent to $v_j$ and $a_{ij}=0$ otherwise. If $G$ and $G'$
represent the same graph-link then
$\operatorname{corank}_{\mathbb{Z}_2}(A(G)+E)=\operatorname{corank}_{\mathbb{Z}_2}(A(G')+E)$,
where $E$ is identity matrix.

 \begin {definition}\label {def:num_com}
Let us define the {\em number of components} in a graph-link
$\mathfrak{F}$ as $\operatorname{corank}_{\mathbb{Z}_2}(A(G)+E)+1$,
here $G$ is a representative of $\mathfrak{F}$. A graph-link
$\mathfrak{F}$ with $\operatorname{corank}_{\mathbb{Z}_2}(A(G)+E)=0$
for any representative $G$ of $\mathfrak{F}$ is called a {\em
graph-knot}.
 \end {definition}

Let $\operatorname{corank}_{\mathbb{Z}_2}(A(G)+E)=0$,
$B_i(G)=A(G)+E+E_{ii}$ (all elements of $E_{ii}$ except for the one
in the $i$-th column and $i$-th row which is one are $0$) for each
vertex $v_i\in V(G)$.

 \begin {definition}\label {def:wr_num}
Let us define the {\em writhe number} $w_i$ of $G$ (with
$\operatorname{corank}_{\mathbb{Z}_2}(A(G)+E)=0$) at $v_i$ as
$w_i=(-1)^{\operatorname{corank}_{\mathbb{Z}_2}
B_i(G)}\operatorname{sign}v_i$ and the {\em writhe number} of $G$ as
  $$
w(G)=\sum\limits_{i=1}^nw_i.
  $$
 \end {definition}

If $G$ is a realisable graph by a chord diagram and, therefore, by a
virtual diagram then $w_i$ is the ``real'' writhe number of the
crossing corresponding to $v_i$.

 \begin {definition}
We say that an $n\times n$ matrix $A=(a_{ij})$ {\em coincides with
an $n\times n$ matrix $B=(b_{kl})$ up to diagonal elements} if
$a_{ij}=b_{ij}$, $i\ne j$.
 \end {definition}

  \begin {lemma}[\cite {Ily1}]\label {lem:nondeg}
Let $A$ be a symmetric matrix over $\mathbb{Z}_2$. Then there exists
a symmetric matrix $\widetilde{A}$ with $\det\widetilde{A}=1$ equal
to $A$ up to diagonal elements.
 \end {lemma}

Let $\mathfrak{F}$ be a graph-knot and let $G$ be its
representative. Let us consider the simple graph $H$ having the
adjacency matrix coinciding with $(A(G)+E)^{-1}$ up to diagonal
elements and construct the graph $L(G)$ from $H$ by just adding
loops to any vertex of $H$ corresponding to a vertex of $G$ with the
negative writhe number. Let us define the map $\chi$ from the set of
graph-knots to the set of homotopy classes of looped interlacement
graphs defined by $\chi(\mathfrak{F})=\mathfrak{L}$, here
$\mathfrak{L}$ is the homotopy class of $L(G)$. It turns out that
the map $\chi$ is well-defined~\cite {Ily2}. The inverse map is
defined as follows~\cite {Ily2}. Let $\mathfrak{L}$ be the homotopy
class of $L$. By using Lemma~\ref {lem:nondeg} we can construct a
symmetric matrix $A=(a_{ij})$ over $\mathbb{Z}_2$ coinciding with
the adjacency matrix of $L$ up to diagonal elements and $\det A=1$.
Let $G(L)$ be the labeled simple graph having the matrix $A^{-1}+E$
as its adjacency matrix (therefore, the first component of the
vertex label is equal to the corresponding diagonal element of
$A^{-1}+E$), the second component of the label of the vertex with
the number $i$ is $w_i(1-2a_{ii})$, here $w_i=1$ if the vertex of
$L$ with the number $i$ doesn't have a loop, and $w_i=-1$
otherwise~\cite {Ily1,Ily2}. We have
$\chi^{-1}(\mathfrak{L})=\mathfrak{F}$, here $G(L)$ is a
representative of $\mathfrak{F}$. Therefore, we get

 \begin {theorem}[\cite {Ily1,Ily2}]
There is a one-to-one correspondence between the set of graph-knots
and the set of homotopy classes of looped interlacement graphs.
Moreover, the graph-knot and the homotopy class of looped
interlacement graphs both constructed from a given virtual knot
diagram are related by this map.
 \end {theorem}

\subsection{Soboleva's Theorem and its Corollaries}\label{sob}

\begin {figure}
\centering\includegraphics[width=300pt]{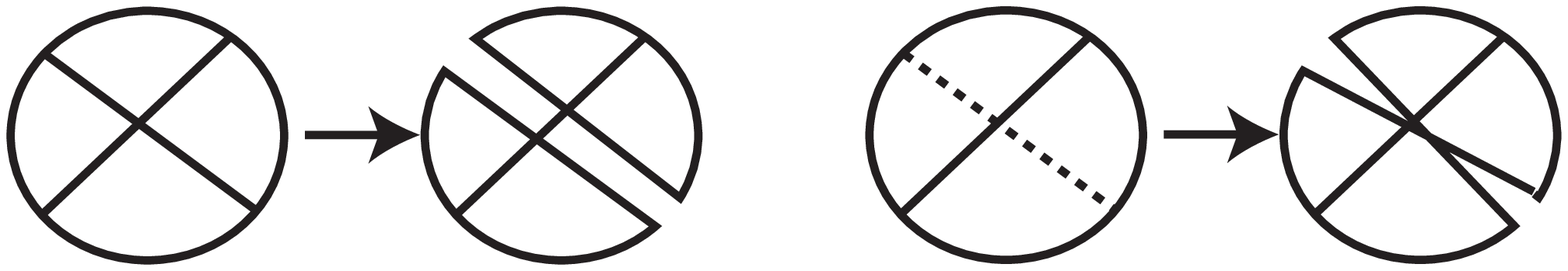}\caption{A surgery
of the circuit along a chord} \label{surger}
 \end {figure}

\begin {figure}
\centering\includegraphics[width=200pt]{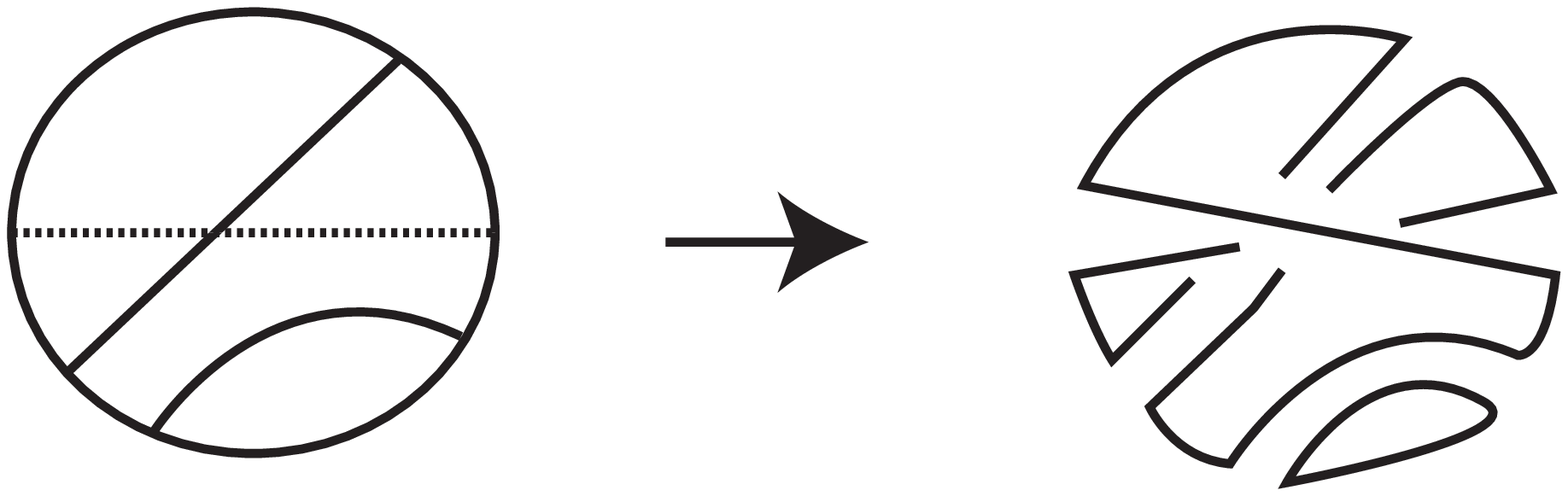}\caption{The
manifold $m(D)$} \label{result}
 \end {figure}

Assume we are given a framed chord diagram, i.e.\ each chord of it
is endowed only with framing. Define the {\em surgery over the set
of chords as follows}. For every chord having the framing $0$
(resp., $1$), we draw a parallel (resp., intersecting) chord near it
and remove the arc of the circle between adjacent ends of the chords
as in Fig.~\ref {surger}. By a small perturbation, the picture in
$\mathbb{R}^{2}$ is transformed into a one-manifold in
$\mathbb{R}^{3}$. This manifold $m(D)$ is the result of surgery, see
Fig.~\ref {result}. Surprisingly, the number of the connected
components of $m(D)$ can be determined from the intersection graph.

  \begin {theorem}[\cite {Soboleva,Tr}]\label {th:sob}
Let $D$ be a labeled chord diagram, and let $G$ be its labeled
intersection graph. Then the number of connected components of
$m(D)$ equals $\operatorname{corank}_{\mathbb{Z}_2}A(G)+1$, where
$A(G)$ is the adjacency matrix of $G$.
  \end {theorem}

 \begin {remark}
Theorem~\ref {th:sob} allows us to define `the number of circles'
after a ``surgery'' for a graph even when the given graph is not an
intersection graph of any chord diagram.

We refer to this Theorem as Soboleva's theorem although a partial
case of this Theorem, when all the chords have the framing $0$,
firstly appears in~\cite {CL}.
 \end {remark}

Throughout the rest of the  paper, we shall use the following {\em
main principle}:

{\em Assume there is an equality concerning numbers of circles in
some states of chord diagrams, which can be formulated in terms of
the intersection graph. Then the corresponding equality usually
holds even for non-realisable graphs.}

The reason for this principle to hold is the following: every time
we have a picture where some two numbers of circles are equal to
each other (or differ by a constant), this can be expressed in terms
of the corresponding adjacency matrices, and the proof does not
generally depend on the behaviour of the matrix outside of the
crossings in question. This means that the equality holds true for
generic matrices, thus, it works for general intersection graphs.

This principle has lots of consequences. We shall demonstrate it for
three examples.

The first example comes from Definition~2.14. It shows that the
number of components defined for non-realisable graphs by the same
formula as for realisable ones doesn't change under the Reidemeister
moves.

The second example comes from Definition~2.15. We have defined the
writhe number of a vertex for an arbitrary labeled graph by using
the definition of the writhe number for virtual diagrams and writing
this definition in terms of matrices. Then we define the writhe
number for a graph. Since the proof exists in the realisable case,
it can be rewritten in terms in matrices, thus, we get that the
writhe number of a graph doesn't change under the second and third
Reidemeister graph-moves and changes by $\pm1$ under the first
Reidemeister graph-moves.

The third example is as follows. Assume we have a framed $4$-graph
$K$ with a vertex $v$ and we would like to know whether this vertex
belongs to one component or it belongs to different components of
the corresponding graph-link. Then we may take the two smoothings
$K_{a},K_{b}$ of the $K$ at $v$ and see how many components we get.
If $v$ belongs to two branches of the same component of $K$ then the
number of components of one of $K_{a},K_{b}$ is equal to that of
$K$, and the number of components of the other one is equal to that
of $K$ plus one. If $v$ belongs to two different components of $K$,
then the number of components of each of $K_{a},K_{b}$ is that of
$K$ minus one. Now, turning to graph-links, by taking appropriate
matrix ranks, we may see whether each vertex belongs to one
component or to two different components of the graph-link. This
method is used for proving the invariance of the Kauffman bracket
polynomial.

\subsection{Smoothing Operations and Turaev's $\Delta$}
\label{smsec}

\begin{figure}
\centering\includegraphics[width=200pt]{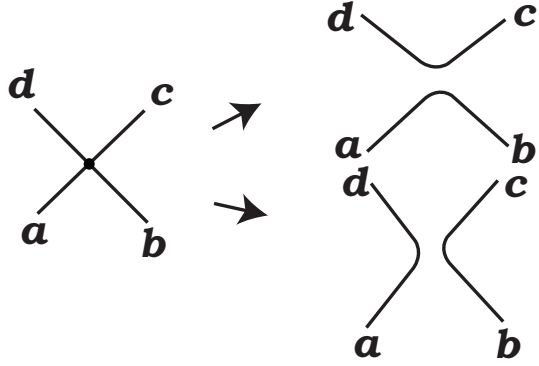}\caption{Two
smoothings at a vertex for a framed graph} \label{smooth12}
\end{figure}

By a {\em smoothing} of a framed $4$-graph at a vertex $v$ we mean
any of the two framed $4$-graphs obtained by removing $v$ and
repasting the edges as $a-b$, $c-d$ or as $a-d,b-c$, see
Fig.~\ref{smooth12}. Generally, a {\em smoothing} of a framed
$4$-graph in a collection of vertices is the framed $4$-graph
obtained by a sequence of smoothings.

\begin{figure}
\centering\includegraphics[width=200pt]{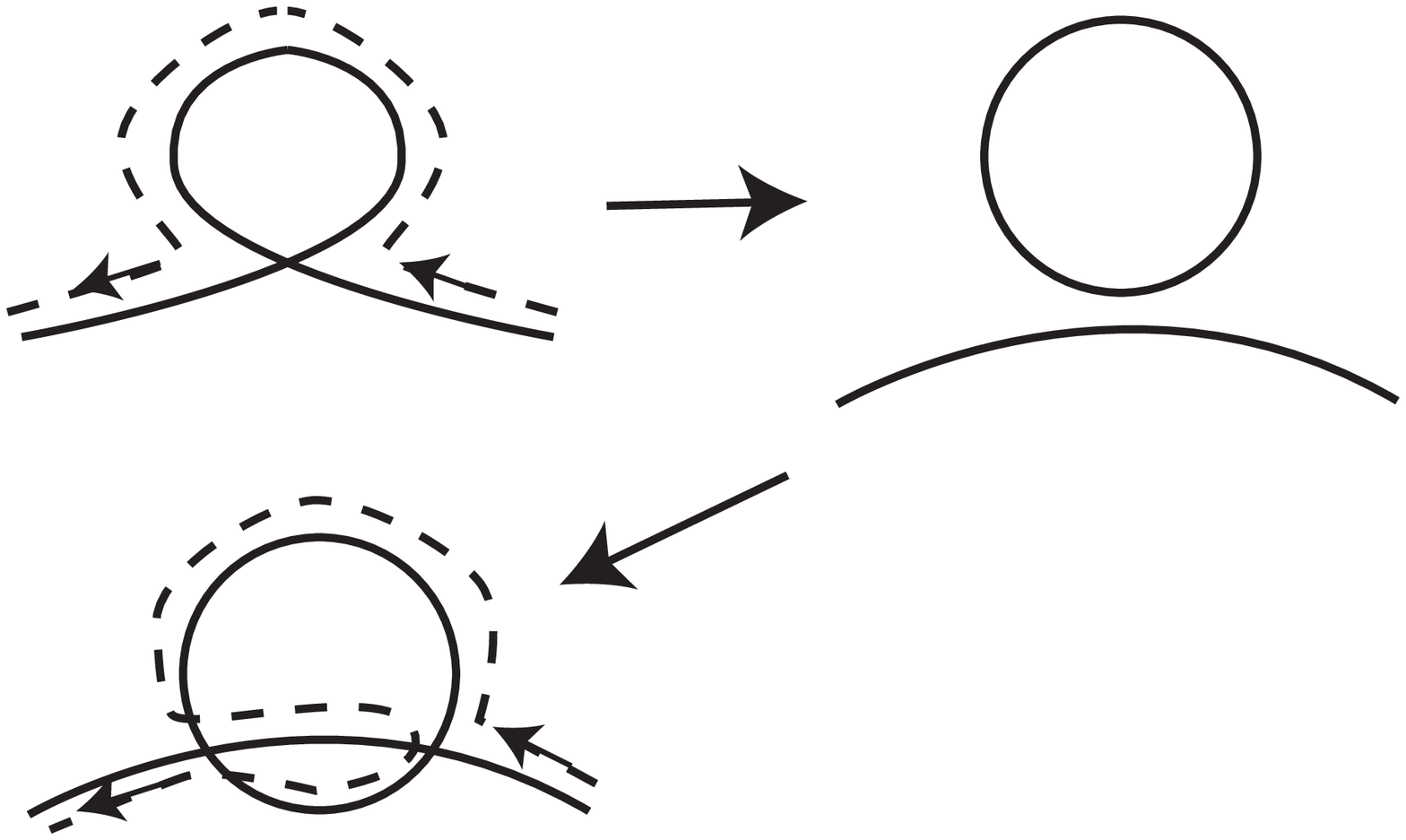} \caption{One of
the two smoothings at an isolated vertex} \label{sm_gr}
\end{figure}

Now, we can mimic this definition for the case of free graph-links.
Let $\mathcal{G}$ be a free framed graph, i.e.\ an equivalence class
of labeled simple graphs, and let $v\in V(\mathcal{G})$ (the set of
vertices is the same for any representative of $\mathcal{G}$). Let
us consider two cases. In the first case there exists a
representative $H$ of $\mathcal{G}$ for which $v$ has either framing
$1$ or the degree more than $0$. It is not difficult to see that $v$
has the same property for each representative of $\mathcal{G}$, and
there are two representatives $H_1$ and $H_2$ of $\mathcal{G}$ which
distinguish from each other by $\Omega_g 4$ or $\Omega_g 4'$ at $v$.
By a {\em smoothing} of a free framed graph $\mathcal{G}$ at a
vertex $v$ we mean any of the two free framed graphs having the
representatives $H_1\setminus\{v\}$ and $H_2\setminus\{v\}$,
respectively. In the second case $v$ has framing $0$ and is isolated
for each representative of $\mathcal{G}$. Let $H$ be a
representative of $\mathcal{G}$. Let us construct the new graph $H'$
obtained from $H$ by adding a new vertex $u$ with framing $0$ to $H$
which is adjacent only to $v$, see Fig.~\ref {sm_gr} for the case of
realisable graphs (the dashed line is a rotating circuit). By a {\em
smoothing} of a free framed graph $\mathcal{G}$ at a vertex $v$ we
mean any of the two free framed graphs having the representatives
$H\setminus\{v\}$ and $H'$, respectively. It is not difficult to see
that the numbers of the components of two different smoothings at
any vertex distinguish from each other by $1$. Generally, a {\em
smoothing} of a free framed graph in a collection of vertices is the
free framed graph obtained by a sequence of smoothings.

We call a free framed graph $\mathcal{G}$ {\em realisable} if any of
its representative is realisable by a chord diagram. It is not
difficult to show that $\mathcal{G}$ is realisable if and only if
there exists a realisable representative of it, we just redraw the
picture.

 \begin{statement}
Assume $\mathcal{G},\,\mathcal{G}'$ are free framed graphs and
$\mathcal{G}$ can be obtained from $\mathcal{G}'$ as a result of
smoothings at some vertices. Then if $\mathcal{G}'$ is realisable by
a chord diagram, then so is $\mathcal{G}$.\label{smreal}
 \end{statement}

The proof is obvious: if $\mathcal{G}'$ is realisable and
$\mathcal{G}$ is obtained from $\mathcal{G}'$ by applying a fourth
graph-move and/or deleting a vertex, then one can draw the
corresponding framed $4$-graph, and take the corresponding
resmoothing which will yield the chord diagram for $\mathcal{G}$. In
the other case, if $\mathcal{G}'$ is realisable then the
realisability of $\mathcal{G}$ follows from Lemma~\ref {lem:co_re}.

Let $i>1$ be a natural number. Define the set
$\mathbb{Z}_2\mathcal{G}_i$ to be the $\mathbb{Z}_2$-linear space
generated by the set of free framed graphs $\mathcal{G}$ with
$\operatorname{corank}_{\mathbb{Z}_2}(A(\mathcal{G})+E)=i-1$ modulo
the following relations:

1) the second Reidemeister graph-moves;

2) $\mathcal{G}=0$, if $\mathcal{G}$ has two vertices with framing
$0$ which are adjacent only to each other.

For $i=1$, we define $\mathbb{Z}_2\mathcal{G}_1$ analogously with
respect to equivalence 1) and not 2).

Let us define the map
$\Delta\colon\mathbb{Z}_2\mathcal{G}_1\to\mathbb{Z}_2\mathcal{G}_1$,
cf.~\cite {Ma3}. Given a free framed graph $\mathcal{G}$ with
$\operatorname{corank}_{\mathbb{Z}_2}(A(\mathcal{G})+E)=0$. We shall
construct an element $\Delta(\mathcal{G})$ from
$\mathbb{Z}_2\mathcal{G}_2$ as follows. For each vertex $v$ of
$\mathcal{G}$, there are two ways of smoothing it. One way gives a
graph from $\mathbb{Z}_2\mathcal{G}_1$, and the other smoothing
gives a free framed graph $\mathcal{G}_v$ from
$\mathbb{Z}_2\mathcal{G}_2$. We take $\mathcal{G}_v$ and set
 $$
\Delta(\mathcal{G}) =\sum\limits_v
\mathcal{G}_v\in\mathbb{Z}_2\mathcal{G}_2.
 $$

 \begin {theorem}
$\Delta$ is a well defined mapping from $\mathbb{Z}_2\mathcal{G}_1$
to $\mathbb{Z}_2\mathcal{G}_2$.
 \end {theorem}

By using the main principle from subsection~\ref {sob} we can define
whether a vertex belongs to one ``component'' or different
components of a free framed graph. Namely, we call a vertex $v_i$ of
a free framed graph $\mathcal{G}$ {\em oriented} if
$\operatorname{corank}_{\mathbb{Z}_2}(A(\mathcal{G})+E)\leqslant\operatorname{corank}_{\mathbb{Z}_2}(B_i(\mathcal{G}))$.
It is not difficult to show that
$\operatorname{corank}_{\mathbb{Z}_2}B_i(\mathcal{G})\ne\operatorname{corank}_{\mathbb{Z}_2}B(\mathcal{G}\setminus\{v_i\})$
if $v_i$ is oriented.

By using the notion of an oriented vertex we can define the map
$\Delta^{i}$ (iteration) by considering smoothings at oriented
vertices and taking that smoothing which has more components than
other in each step.

 \begin {corollary}
$\Delta^{i}$ is a well defined mapping from
$\mathbb{Z}_2\mathcal{G}_1$ to $\mathbb{Z}_2\mathcal{G}_{i+1}$.
 \end {corollary}

\section{Parity, Minimality and Non-Triviality Examples}

In the present section we consider the parity for free graph-knots
and free graph-links in spirit of~\cite{Ma1,Ma2}. As we have
constructed the one-to-one correspondence between the set of
graph-knots and the set of homotopy classes of looped interlacement
graphs it is sufficient to construct a parity for free homotopy
classes of looped interlacement graphs and for graph-links with more
than one components.

Consider the category of free homotopy classes of looped
interlacement graphs.

 \begin {definition}
For every looped interlacement graph $L$ we call a vertex $v$ {\em
even} if $v$ is adjacent to an even number of vertices distinct from
$v$; otherwise we call $v$ {\em odd}.
 \end {definition}

 \begin {theorem}
The parity defined on homotopy classes of looped interlacement
graphs {\em (}resp., free homotopy classes\/{\em)} by using even and
odd vertices satisfies the parity axioms from~\cite
{Ma3}.\label{parityknots}
 \end {theorem}

The proof evidently follows from the consideration of the
Reidemeister moves for looped interlacement graphs.

By using this parity we can define the map $\Delta^i_{odd}$ where
the sum is taken over all odd oriented vertices or $\Delta^i_{even}$
where the sum is taken over all even oriented vertices. We have to
define the notion of {\em even} and {\em odd} vertex for free framed
graphs with many components.

We call a vertex $v$ of $\mathcal{G}$ with one component {\em even}
(resp., {\em odd}) if the vertex corresponding to $v$ of the looped
interlacement graph $\chi(\mathcal{G})$ is even (resp., odd). Let us
consider the free framed graph $\mathcal{G}_{v_1,\dots,v_{k-1}}$
with $k$ components which is obtained from $\mathcal{G}$ by
smoothing $\mathcal{G}$ consequently at $v_1,\dots,v_{k-1}$; $v_i$
is a oriented vertex in $\mathcal{G}_{v_1,\dots,v_{i-1}}$. An
oriented vertex $u$ of $\mathcal{G}_{v_1,\dots,v_{k-1}}$ is {\em
even with respect to the smoothing at $v_1,\dots,v_{k-1}$} (resp.,
{\em odd with respect to the smoothing at $v_1,\dots,v_{k-1}$}) if
the number of oriented vertices in $\mathcal{G}_{v_1,\dots,v_{k-1}}$
which are incident to $u$ in $\chi(\mathcal{G})$ is even (resp.,
odd).

 \begin {remark}
We have defined even vertices only for those free framed graphs with
many components which originate from given free framed graphs with
$1$ component. It is sufficient to define the iteration
$\Delta^i_{odd}$ and $\Delta^i_{even}$.
 \end {remark}

 \begin {statement}
$\Delta^{i}_{odd}$ is a well defined mapping from
$\mathbb{Z}_2\mathcal{G}_1$ to $\mathbb{Z}_2\mathcal{G}_{i+1}$.
 \end {statement}

Let us consider another parity for graph-links with two components.

 \begin {definition}
For a graph $G$ representing a two-component graph-link we call a
vertex {\em even} if it is {\em oriented} and {\em odd} otherwise.
 \end {definition}

 \begin {theorem}
The parity defined above satisfies the parity axioms given in~\cite
{Ma3}.\label{paritylinks}
 \end {theorem}

The proof again follows from the definition of parity.

 \begin{remark}
The parity used in Theorem~\ref{parityknots} is easier to define via
{\em Gauss diagram approach}: assuming a graph is the ``intersection
graph of a non-existing Gauss diagram'', we take those chords which
are linked with even number of chords to be even, and the remaining
ones to be odd.

In the language of the rotating circuit approach, this is more
difficult: just taking an arbitrary rotating circuit and counting
the number of adjacent chords can lead to different results: the
same vertex may have an even or odd degree for different rotating
circuits.
\end{remark}

Now, having these parities in hand, we can define the brackets
$[\cdot]$ for graph-knots and $\{\cdot\}$ for graph-links
analogously to the case of~\cite {Ma3}.  For a free graph-knot $G$
(resp., free two-component graph-link $H$), consider the following
sums
 $$
[G]=\sum_{s\;even.,1\; comp} G_{s},
 $$
and
 $$
\{H\}=\sum_{s\;even.\; non-trivial} H_{s},
 $$
where the sums are taken over all smoothings at all {\em even}
vertices (with respect to the definitions above), and only those
summands are taken into account where
$\operatorname{corank}_{\mathbb{Z}_2}(A(G_{s})+E)=0$ (resp., $H_{s}$
is not equivalent to any simple graph having two vertices with
framing $0$ which are adjacent only to each other). Thus, if $G$ has
$k$ even vertices, then $[G]$ will contain at most $2^{k}$ summands,
and if all vertices of $G$ are odd, then we shall have exactly one
summand, the graph $G$ itself. The same is true for $H$ and $\{H\}$.

Now, we are ready to formulate the main theorems of this section:

  \begin{theorem}
If $G$ and $G'$ represent the same free graph-knot then the
following equality holds: $[G]=[G']$.

Analogously, if $H$ and $H'$ represent the same free graph-link with
two components then $\{H\}=\{H'\}$. \label {mainthm}
  \end{theorem}

The proof of~\ref {mainthm} verbally reproduces the proof of the
main theorem from~\cite {Ma3}, according to the {\em main principle}
or, maybe, a slight modification of it.

 \begin {definition}
We call a labeled graph $G$ (resp., a looped graph $L$) {\em
minimal} if there is no representative of the graph-link
corresponding to $G$ (resp., the homotopy class of $L$) having
strictly smaller number of vertices than $G$ (resp., $L$) has.
 \end {definition}

  \begin{theorem}\label {th:main_pair}
Let $G$ {\em (}resp., $H$\/{\em)} be a simple labeled graph
representing a free graph-knot {\em(}a two-component
graph-link\/{\em)} with all odd vertices in the sense of
Theorem~\ref {parityknots} {\em(}resp., Theorem~\ref
{paritylinks}{\em)}, such that no decreasing second Reidemeister
move is applicable to $G$ {\em(}resp., $H${\em)}. Then there is no
simple graph equivalent to $G$ {\em(}resp., $H${\em)} with strictly
smaller number of vertices.\label{brackets}
  \end{theorem}

As a consequence of this theorem we may deduce the following

  \begin{corollary}
The free graph-knot $G$ shown in Fig.~\ref{Bouchet9} is minimal; in
particular, it is non-trivial and has no realisable representatives.

Moreover, $G$ has no representative realisable as an intersection
graph of a chord diagram.\label{maincrl}
  \end{corollary}

\begin{figure}
\centering\includegraphics[width=200pt]{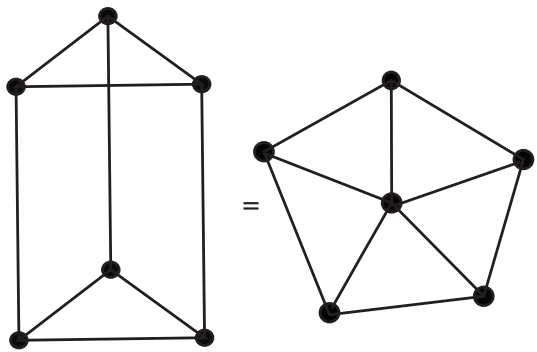} \caption{The
Two-Component Free Link} \label {prism13}
 \end{figure}

The first claim of the corollary is trivial: we just check the
conditions of the theorem and see that $G$ is minimal; to see that
the second claim indeed holds, we shall need to go through the proof
of Theorem~\ref {mainthm}, where we see that any representative $G'$
of the graph-knot $\mathbb{G}$ has $G$ as a {\em smoothing}, that
is, $G$ {\em lies inside} each representative $G'$ of the same
graph-link, and if $G$ is not realisable, then so is $G'$ according
to Statement~\ref {smreal}. Analogously, one sees that the free
two-component graph-link $\mathfrak{H}$ with the representative $H$
shown in Fig.~\ref {prism13} (left part) has no realisable
representative because $\{H\}=H$. Note that $H$ in Fig.~\ref
{prism13} is equivalent to the Bouchet graph shown in the right part
of the picture by $\Omega_{4}$-graph move; so they represent the
same two-component free graph-link.

Let us consider one more example.

  \begin{statement}
The looped graph $K$ represented by the ``Gauss diagram'' shown in
Fig.~\ref {xx} is minimal and non-realisable.
  \end{statement}

The proof consists of the following steps very similar to the
example from~\cite {Ma3}.

First, note that $\Delta(K)$ consists of $7$ summands
$L+\sum_{i}L_{i}$, where only one summand (corresponding to the
vertex $x$) is a $2$-component free graph-link with all {\em odd}
vertices; for each of the remaining summands $L_i$, there is at
least one even vertex. This is, indeed, very easy: the graph $K$ has
only one vertex $x$ which is adjacent to all the remaining vertices,
then the argument just repeats the argument from~\cite {Ma3}.

\begin{figure}
\centering\includegraphics[width=300pt]{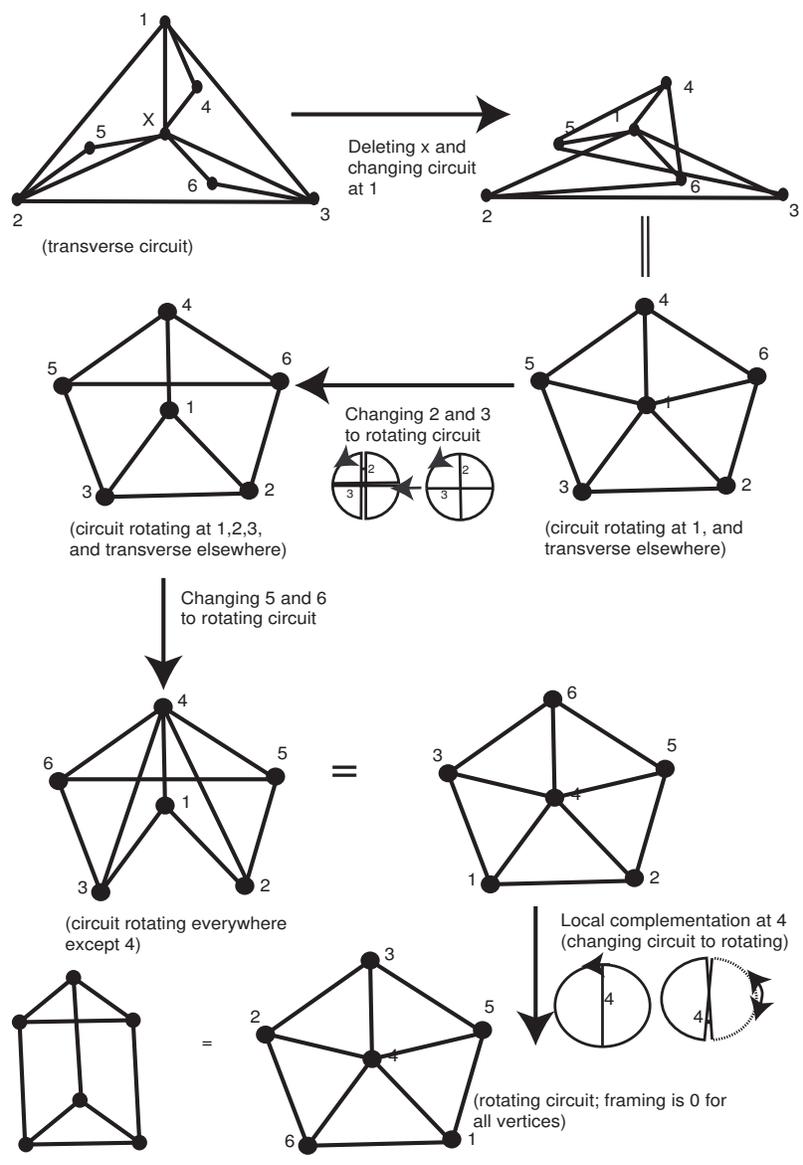} \caption{The
Two-Component Free Link} \label {ktoprism14}
 \end{figure}

Now, the $2$-component free graph-link $L$ has a rotating circuit
diagram shown in Fig.~\ref {prism13}; all framings of the vertices
are $0$. To see it, one should consecutively perform the following
operations for $K$: first, we ``smooth'' its crossing $x$; it can be
done only at the expense of changing our circuit to the rotating one
at some vertex different from $x$ (a Gauss circuit can not represent
link); after that, we have to make the circuit rotating at all other
vertices. All these steps are shown in Fig.~\ref {ktoprism14}.

Now, consider the bracket $\{\Delta(K)\}=L+\sum_{i}\{L_{i}\}$. Note
that all summands $\{L_{i}\}$ have representatives with strictly
less than $6$ vertices since each of $L_{i}$ has at least one even
vertex; on the other hand, $L$ has no representative with less than
$6$ crossings; so, this element $L$ is not canceled  in the sum.
Since it is not realisable, the free framed knot $K$ is not
realisable either.

\section{Atoms and Orientability}

In this section we introduce a object: {\em atom} which is a
$2$-manifold with additional structure. This object helps us define
new equivalence relation of graphs, see~\cite {IM1}, and is very
useful for minimality theorems, see ahead.

  \begin {definition}
An {\em atom}~\cite {F} is a pair $(M,\Gamma)$ consisting of a
closed $2$-manifold $M$ and a finite graph $\Gamma$ embedded in $M$
together with a colouring of $M\backslash \Gamma$ in a checkerboard
manner. An atom is called {\em orientable} ({\em connected}) if the
surface $M$ is orientable (connected). Here $\Gamma$ is called the
{\em frame} of the atom. By {\em genus} ({\em Euler characteristic,
orientation}) of the atom we mean that of the surface $M$. Atoms and
their genera were also studied by Turaev~\cite{Turg}, and in
consequent papers that atom genus is also called the Turaev
genus~\cite {Turg}.
  \end {definition}

Having a chord diagram, we can construct an atom corresponding to
the chord diagram, see~\cite{IM1,IM2}.

In fact, chord diagrams in the sense of rotating circuits with all
chords having framing $0$ encode all orientable atoms, see~\cite
{Ma4}. Chord diagrams with all positive chords encode all atoms with
one white cell: this white cell corresponds to the $A$-state of the
virtual diagram, and chords show how this cell approaches itself in
neighbourhoods of crossings (atom vertices). If we want to deal with
all atoms and restrict ourselves for the case of one circle, we
should take this circle to correspond to some other state of the
atom, which is encoded by labelings of the chords.

\section{Kauffman's Bracket Generalisation and Other Invariants. Minimality Theorems}

We have already considered some minimality theorems. In this section
we present minimality theorems which use Kauffman's bracket
generalisation. Establishing minimal crossing number of a certain
link is one of the important problems in classical knot theory. In
late 19's century, famous physicist and knot tabulator
P.\,G.~Tait~\cite {Tait} conjectured that alternating prime diagrams
of classical links are minimal with respect to the number of
classical crossings. This celebrated conjecture was solved only in
1987, after the notions of the Jones polynomial and the Kauffman
bracket polynomial appeared. The first solution was obtained by
Murasugi~\cite{Mur}, then it was reproved by
Thistlethwaite~\cite{Thi1}, Turaev~\cite{Turg}, and others. Later,
Thistlethwaite~\cite{Thi2} established the minimality for a larger
class of diagrams (so-called {\em adequate} diagrams). It turns out
that many results in these directions generalise for virtual links
(establishing the minimal number of classical crossings); these
results were obtained by the second named author of the present
paper. See~\cite{MBook} for the proofs and some further
generalisations and other results concerning virtual knots.

 \begin {definition}
The difference between the leading degree and the lowest degree of
non-zero terms of a polynomial $P(x)$ is called the {\em span} of
$P(x)$ and is denoted by $\operatorname{span} P(x)$.
 \end {definition}

 \begin {definition}
A classical link diagram is called {\em alternating} if while
passing along every component of it we alternate undercrossings and
overcrossings.
 \end {definition}

From the `atomic' point of view, alternating link diagrams are those
having atom genus (Turaev genus) zero (more precisely, diagram has
genus zero if it is a connected sum of several alternating
diagrams).

 \begin {definition}
A virtual link diagram $D$ is called {\em split} if there is a
vertex $X$ of the corresponding atom $(M,\Gamma)$ such that
$\Gamma\backslash X$ is disconnected.
 \end {definition}

The main reason of these minimality theorems and further crossing
estimates come from the well-known Kauffman-Murasugi-Thistlethwaite
Theorem:

 \begin {theorem}
For a non-split classical link diagram $K$ on $n$ crossings we have
$\operatorname{span}\langle K\rangle\leqslant 4n$, whence for
alternating non-split diagrams we have $\operatorname{span}\langle K
\rangle=4n$. Here, $\langle K\rangle$ is the Kauffman bracket of
$K$.
 \end {theorem}

Note that the span of the Kauffman bracket is invariant under all
Reidemeister moves. This theorem is generalised for virtual
diagrams~\cite{MBook}. The estimate $\operatorname{span}\langle
K\rangle\leqslant 4n$ can be sharpened to
$\operatorname{span}\langle K\rangle \leqslant 4n-4g$, where $g$ is
the genus of the corresponding atom.

All the minimality theorems for the case of graph-links in this
section rely on the generalisation of the Kauffman bracket
polynomial for graph-links. Now, let us generalise the notions
defined above for the case of graph-links and define the Kauffman
bracket polynomial for a labeled graph $G$.

Let $s\subset V(G)$ be a subset of the set $V(G)$ of vertices of
$G$. Set $G(s)$ to be the induced subgraph of the graph $G$ with the
set of vertices $V(G(s))=s$ and the set of edges $E(G(s))$ such that
$\{u,v\}\in E(G(s))$, where $u,v\in s$, if and only if $\{u,v\}\in
E(G)$.

 \begin {definition}
We call a subset of $V(G)$ a {\em state} of the graph $G$. The
$A$-{\em state} is the state consisting of all the vertices of $G$
labeled $(a,-)$, $a\in\{0,1\}$, and no vertex labeled $(b,+)$,
$b\in\{0,1\}$. Analogously, the $B$-{\em state} is the state
consisting of all vertices of $G$ labeled $(b,+)$, $b\in\{0,1\}$,
and no vertex labeled $(a,-)$, $a\in\{0,1\}$.
 \end {definition}

 \begin {definition}
The {\em Kauffman bracket polynomial} of $G$ is
  $$
\langle
G\rangle(a)=\sum\limits_{s}a^{\alpha(s)-\beta(s)}(-a^2-a^{-2})^{\operatorname{corank}
A(G(s))},
  $$
where the sum is taken over all states $s$ of the graph $G$,
$\alpha(s)$ is equal to the sum of the vertices labeled $(a,-)$,
$a\in\{0,1\}$, from $s$ and the vertices labeled $(b,+)$,
$b\in\{0,1\}$, from $V(G)\setminus s$, $\beta(s)=|V(G)|-\alpha(s)$.
 \end {definition}

 \begin {theorem}\label {th:mov_reid}
The Kauffman bracket polynomial of a labeled graph is invariant
under $\Omega_g 2 - \Omega_g 4'$ graph-moves and is multiplied by
$(-a^{\pm3})$ under $\Omega_g 1$ graph-move.
 \end {theorem}

 \begin {definition}
For a labeled graph $G$ let the {\em atom genus} ({\em Turaev
genus}) be $1-(k+l-n)/2$, where $k$ and $l$ are the numbers of
circles in the $A$-state $s_1$ and the $B$-state $s_2$ of $G$,
respectively, i.e.\ $k=\operatorname{corank}_{\mathbb{Z}_2}
A(G(s_1))+1$ and $l=\operatorname{corank}_{\mathbb{Z}_2}
A(G(s_2))+1$.
 \end {definition}

Note that this number agrees with the atom genus in the usual case:
we just use $\chi=2-2g$, where $\chi$ is the Euler characteristic,
and count $\chi$ by using the number of crossings $n$, number of
edges $2n$ and the number of $2$-cells ($A$-state circles and
$B$-state circles).

 \begin {definition}
A labeled graph $G$ on $n$ vertices is {\em alternating} if its atom
genus is equal to $0$. A labeled graph $G$ is {\em non-split} if it
has no isolated vertices.
 \end {definition}

 \begin {theorem}\label {th:minimal}
An alternating non-split labeled graph is minimal.
 \end {theorem}

\begin{figure}
\centering\includegraphics[width=200pt]{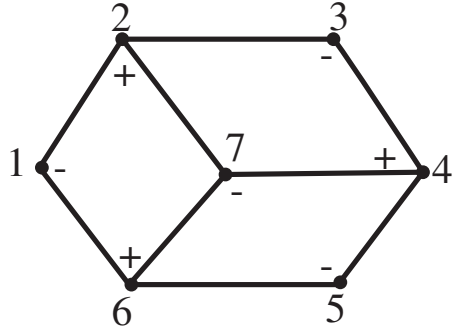}\caption{The graph
$BW_3$} \label{fig:g7}
 \end{figure}

{\bf Example.} Consider the graph $BW_{3}$ consisting of the $7$
vertices with the following incidences. For $i,j=1,\dots, 6$, $i$ is
connected to $j$ if and only if $i-j\equiv\pm1\pmod 6$, and $7$ is
connected to $2,4,6$. Label all even vertices by $(0,+)$, and label
all odd vertices by $(0,-)$, see Fig.~\ref {fig:g7}. This graph is
alternating. By Theorem~\ref {th:minimal}, $BW_{3}$ is minimal. Note
that this graph is not realisable as an intersection graph of a
chord diagram. We conjecture that the graph-link represented by
$BW_{3}$ is non-realisable. At least we know that it has no such
representatives with the number of crossings less than or equal to
$7$.

We conclude the present paper with the significant one question
which seems to be very interesting for us: is there a graph-link
having two representatives realisable by chord diagrams but they are
equivalent in the graph-link only by means of non-realisable graphs?

\end{document}